

\baselineskip=14pt
\parskip=10pt

\font\eightrm=cmr8 
\font\eighttt=cmtt8
\magnification=\magstephalf

\def\1{{\overline{1}}}
\def\2{{\overline{2}}}
\parindent=0pt
\overfullrule=0in

\def\frac#1#2{{#1 \over #2}}
\bf
\centerline
{
On Euler's ``Misleading Induction", Andrews' ``Fix", and How to Fully Automate them
}
\rm
\bigskip
\centerline{ {\it Shalosh B. EKHAD and
Doron 
ZEILBERGER}\footnote{$^1$}
{\eightrm  \raggedright
Department of Mathematics, Rutgers University (New Brunswick),
Hill Center-Busch Campus, 110 Frelinghuysen Rd., Piscataway,
NJ 08854-8019, USA.
{\eighttt zeilberg  at math dot rutgers dot edu} ,
\hfill \break
{\eighttt http://www.math.rutgers.edu/\~{}zeilberg/} .
April 3, 2013.
Accompanied by Maple package \hfill \break {\eighttt GEA}
downloadable from Zeilberger's website.
Supported in part by the NSF.
}
}

\qquad\qquad\qquad\qquad\qquad\qquad\qquad\qquad\qquad
{\it Dedicated to George Andrews on  his $(75-\epsilon)^{th}$ birthday}

Recall that the {\it trinomial coefficient}([Wei])
$$
{{n} \choose {j}}_2
$$ 
is the coefficient of $x^j$ in 
$$
(1+x+x^{-1})^n \quad .
$$
In other words, 
$$
(1+x+x^{-1})^n=\sum_{j=-n}^{n} {{n} \choose {j}}_2 x^j \quad .
$$
Also recall that the {\it Fibonacci numbers}, $F_n$, are defined by
$F_{-1}=1, F_0=0$, and $F_n=F_{n-1}+F_{n-2}$ for $n>0$.

The fascinating story of how Euler {\it almost} fooled himself into believing that
$$
3{{n+1} \choose {0}}_2 - {{n+2} \choose {0}}_2 = F_{n}(F_{n}+1) \quad,
\eqno(Leonhard)
$$
for {\bf all} $n$ because he checked this for the {\bf nine} values $-1 \leq n \leq 7$, only to find out
that it fails for $n=8$, leading him to record it for {\it posterity} in [E], has been told several times,
including  the nice `popular' book by David Wells[Wel], Eric Weisstein's  extremely useful Mathworld website[Wei],
and Ed Sandifer's famous on-line MAA column[S].

In the first three sections of George Andrew's important article [An] (that merely serve as
the motivation and background for the remaining sections that talk about deep $q$-analogs), 
he describes a brilliant way to `correct' the left side of of $(Leonhard)$ in order to make the identity come true
for {\it all} $n \geq -1$.

First he used the obvious fact that
$$
{{n+2} \choose {0}}_2={{n+1} \choose {-1}}_2+{{n+1} \choose {0}}_2+{{n+1} \choose {1}}_2
$$
(and symmetry) to rewrite Eq. $(Leonhard)$ as:
$$
{{n+1} \choose {0}}_2 - {{n+1} \choose {1}}_2 = \frac{1}{2}F_{n}(F_{n}+1) \quad,
\eqno(Leonhard')
$$

and then went on to prove (using ad-hoc human ways) the identity
$$
\sum_{j=-\infty}^{\infty} {{n+1} \choose {10j}}_2 - \sum_{j=-\infty}^{\infty} {{n+1} \choose {10j+1}}_2  = \frac{1}{2} F_{n}(F_{n}+1) \quad .
\eqno(George)
$$
Note that for $n<8$ the only non-zero summands in $(George)$ are with $j=0$.

[{\eightrm At the risk of giving away the punch-line,
let's remark that once conjectured, a fully rigorous proof of $(George)$ can be obtained by checking it,
\`a la Euler, for (to be safe) $0 \leq n \leq 20$.}]

{\bf Interlude: Even Giants make stupid conjectures}

We are a little surprised that Euler could have believed, {\it even for a second}, that 
Eq. $(Leonhard)$ is true for {\it all} $n$. Completely by hand (see [S]) Euler
found a three-term linear recurrence with {\it polynomial} coefficients
for ${{n}\choose {0}}_2$ (that easily implies such a recurrence for $3{{n+1}\choose {0}}_2-{{n+2}\choose {0}}_2$,
both are easily found today with the {\it Almkvist-Zeilberger algorithm} [AlZ]),
so he should have realized that it can't equal $F_n(F_n+1)$ that satisfies a linear recurrence
with {\it constant} coefficients (in other words it is what is called today a $C$-finite sequence, see [Z]).

Another way Euler could have easily realized that $(Leonhard)$ is false  is via asymptotics,
even a very crude one. The ratio of consecutive terms on the left side of $(Leonhard)$ obviously tends to $3$ while
the ratios of consecutive terms on the right side tends to $\phi^2=2.61803\dots$.

{\bf The General case}

With Maple (or Sage, or any computer algebra system), it is a piece of cake to generate many Euler-style cautionary tales,
and Andrews-style fixes. Let's summarize our findings by stating a general theorem, whose {\bf proof} also tells you an
{\bf algorithm} how to compute rational generating functions for these sequences. This algorithm has been implemented in
the Maple package {\tt GEA} available from 

{\eighttt http://www.math.rutgers.edu/\~{}zeilberg/tokhniot/GEA} .

{\bf Theorem}: Let $P(x)$ be any Laurent polynomial, and let 
$$
{{n} \choose {j}}_{P}
$$
be the coefficient of $x^j$ in  $P(x)^n$.

Let $k$ be a positive integer, and let $a$ be an integer such that $0\leq a <k$. Define
the {\it generalized Euler-Andrews} sum to be
$$
A(n,k,a):=\sum_{j=-\infty}^{\infty} {{n} \choose {kj+a}}_P \quad .
$$
The {\it generating functions}
$$
f_{k,a}(t):=\sum_{n=0}^{\infty} A(n,k,a)t^n \quad , \quad (0 \leq a <k) \quad ,
$$
are {\bf rational functions} of $t$, all with the {\bf same} denominator, of degree $k$ in $t$.
They are easily computable by linear algebra.

Equivalently, the $k-1$ sequence $\{A(n,k,a)\}_{n=0}^{\infty}$ ($0 \leq a \leq k-1$) satisfy the {\bf same}
homogeneous linear recurrence equation with {\bf constant} coefficients of order $k$ (but of course with
(usually) different initial conditions).

{\bf Proof:} Let's spell-out $P(x)$
$$
P(x)=\sum_{i=\alpha}^{\beta} c_i x^i \quad,
$$
where $\alpha<\beta$ (and $\alpha$ may be negative, of course). Then, obviously, we have the analog of Pascal's triangle identity:
$$
{{n} \choose {j}}_P=
\sum_{i=\alpha}^{\beta} c_i  {{n-1} \choose {j-i}}_P \quad .
$$
Hence
$$
A(n,k,a)=\sum_{i=\alpha}^{\beta} c_i A(n-1,k,a-i \,\,\, mod \,\,\, k)  \quad.
$$
On the level of generating functions we get
$$
f_{k,a}(t):=\delta_{a,0}+t\sum_{i=\alpha}^{\beta} c_i f_{k, (a-i)\,\, mod \,\,\, k}  \quad, \quad (0 \leq a < k) \quad ,
$$
where $\delta_{a,0}$ is $1$ when $a=0$ and $0$ otherwise.
This gives us a system of $k$ linear equations in the $k$ unknowns 
$$
\{ f_{k,0}(t), f_{k,1}(t),  \dots f_{k,k-1}(t) \} \quad,
$$
that Maple can solve very fast.  The fact that the (same) denominator has degree $k$, (and the numerators have
degree $k-1$) follows from Cramer's rule.

This is implemented in procedure {\tt GA(P,x,k,t)}, in the Maple package {\tt GEA}, that
inputs a Laurent polynomial $P$ in the variable $x$, a positive integer $k$ and a variable $t$,
and outputs a list of rational functions in $t$, of length $k$, whose $(a+1)$-th entry is
$f_{k,a}(t)$. If $P$ is symmetric ($P(x)=P(1/x)$) one only has to go as far as $a \leq k/2$,
since then $f_{k,k-a}(t)=f_{k,a}(t)$. This (faster) case is handled by  the procedure {\tt GAs(P,x,k,t)}.

{\bf Computerized Redux of Andrews's man-made proof}

The case $P(x)=x^{-1}+1+x$ and $k=10$ is the one that Andrews needed.
So all you need is type

{\tt GAs(x+1+1/x,x,10,t);} \quad ,

giving all the six generating functions, whose coefficients, $A(m,10,a)$ ($0 \leq a \leq 5$),
are given explicitly in Eq. (2.18) of [An]
(reproduced in [Wei]) as expressions involving  Fibonacci numbers and powers of $3$.
It follows from the theorem (even {\bf without} actually computing the generating functions!)
that Andrews'  claimed formulas can be {\bf proved rigorously} by `just'
checking the first $20$ special cases, as remarked above.

Ditto for Theorem 2.1 of [An], (also quoted in [Wei]). An empirical  proof \`a la Euler (and now
one can manage with $m \leq 10$) suffices.

But if you do not know beforehand  conjectured expressions, then procedure {\tt GA} can find the
generating functions {\it ab initio}.

{\bf The Beauty of Programming}

Even Euler and Andrews would soon get tired of doing the analogous thing for other $k$.
Andrews also did the case $k=6$ in Theorem 3.1 of [An], but we can do it for {\it all} $k$
up to $100$ (easily) and not just for $P(x)=x^{-1}+1+x$, but for {\it any} $P(x)$.
See the sample  output in the front of this article

{\eighttt http://www.math.rutgers.edu/\~{}zeilberg/mamarim/mamarimhtml/gea.html} .

In particular, we can generate ({\bf many!}) Euler-Style {\it cautionary tales} about the premature use of empirical induction, see
procedure {\tt BCT} and the webbook

{\eighttt http://www.math.rutgers.edu/\~{}zeilberg/tokhniot/oGEA6}  \quad .

{\bf Encore: Probabilistic Interpretation}

Suppose you have a (fair or) loaded die whose faces are marked with dollar amounts (some positive, some negative,
some (possibly) $0$),
at each throw you `gain' the amount on the landed face (of course if the amount is $0$ you get nothing,
and if the amount is negative, you have to pay). Let $a(n)$ be the probability of  breaking even after
$n$ throws. Then of course the generating function of $a(n)$ is {\bf not} rational, i.e. the sequence
$a(n)$ does {\bf not} satisfy a linear recurrence with {\bf constant} coefficients
(on the other hand it does satisfy a linear recurrence with {\bf polynomial} coefficients, easily found by
the Almkvist-Zeilberger[AlZ] algorithm, implemented in 
the Maple package {\eighttt http://www.math.rutgers.edu/\~{}zeilberg/tokhniot/EKHAD}.).

But fixing $k$ (even a very large one, say a googol), then the related probability, let's call it $b_k(n)$
of getting an exact multiple of $k$ (that for small $n$ is the same as breaking even),
does satisfy a linear recurrence equation with {\bf constant} coefficients of order $k$,
(equivalently, the generating function is a rational function of degree $k$). 
Ditto for the probability $b_{k,a}(n)$ of finishing with an amount that leaves remainder $a$ when divided by $k$.

{\bf References}

[AlZ] Gert Almkvist and Doron Zeilberger, {\it The Method of Differentiating Under The Integral Sign},
J. Symbolic Computation {\bf 10}(1990), 571-591. \hfill\break
{\tt http://www.math.rutgers.edu/\~{}zeilberg/mamarim/mamarimhtml/duis.html}

[An] George Andrews, {\it Euler's 'exemplum memorabile inductionis fallacis' and q-Trinomial Coefficients},
 J. Amer. Math. Soc. {\bf 3} (1990), 653-669. \hfill\break
{\tt http://www.ams.org/jams/1990-03-03/S0894-0347-1990-1040390-4/S0894-0347-1990-1040390-4.pdf}

[E] Leonhard Euler  {\it Exemplum Memorabile Inductionis Fallacis},
 Opera Omnia, Series Prima, {\bf 15} (1911), 50-69, Teubner. Leipzig, Germany. \hfill\break

[S] Ed Sandifer,
{\it A Memorable Example of False Induction}, MAA on-line \hfill\break
{\tt http://www.maa.org/news/howeulerdidit.html}, 2005.

[Wei] Eric W.  Weisstein, 
{\it Trinomial Coefficient}, From MathWorld--A Wolfram Web Resource. {\tt http://mathworld.wolfram.com/TrinomialCoefficient.html}

[Wel] David Wells, {\it ``Games and Mathematics''}, Cambridge University Press, 2012.

[Z] Doron Zeilberger, {\it The C-finite Ansatz}, to appear in the Ramanujan J. \hfill\break
{\tt http://www.math.rutgers.edu/\~{}zeilberg/mamarim/mamarimhtml/cfinite.html}

\end